\pdfoutput=1
\documentclass[9pt]{article}
\usepackage[utf8]{inputenc}
\usepackage[english]{babel}

\usepackage{eso-pic}
\usepackage{amsmath} 
\usepackage{amssymb}
\usepackage{mathrsfs}
\usepackage{amsthm}
\usepackage{remreset}
\usepackage{tikz}
\usetikzlibrary{calc} 
\usetikzlibrary{math}
\usepackage{tikz-cd} 
\usepackage{graphicx}
\usetikzlibrary{shapes.misc}
\usepackage{appendix}
\usepackage{xcolor}
\usepackage{graphicx}
\usepackage{lineno}
\usepackage{titlesec}
\usepackage{scalerel}[2016/12/29]
\usepackage[numbers]{natbib}

\usepackage{subfiles}
\usepackage[all,cmtip]{xy}

\usepackage{setspace}

\usepackage{hyperref}

\newcommand{\poubelle}[1]{}

\newcommand{\N}{\mathbb{N}}

\newcommand{\U}{\mathcal{U}}

\newcommand{\F}{\mathbb{k}}
\newcommand{\Fp}{\mathbb{F}_p}

\newcommand{\K}{\mathcal{K}}

\newcommand{\E}{\mathcal{V}^f}
\newcommand{\Ef}{\mathcal{V}}

\newcommand{\Hom}{\text{Hom}}

\newcommand{\HomK}{\text{Hom}_\K}
\newcommand{\HomF}{\text{Hom}_{\F}}

\newcommand{\Nil}{\mathcal{N}il}

\newcommand{\id}{\text{id}}

\newcommand{\coker}{\text{coker}}

\newcommand{\SVF}{\mathcal{S}\text{et}^{(\mathcal{V}^f)^\text{op}}}

\newcommand{\IWV}{I_{(W\oplus V,\psi\boxplus \epsilon_V)}}

\newcommand{\FVF}{\mathcal{F}\text{in}^{(\mathcal{V}^f)^\text{op}}}

\newcommand{\decal}{\Bar{\Delta}_{(\F,\epsilon_{\F})}}

\newcommand{\Funcat}{\mathcal{F}(\mathfrak{S}_S,\Ef)}
\newcommand{\Sector}{\mathfrak{S}_S}
\newcommand{\Rector}{\mathfrak{R}_S}
\newcommand{\pol}[1]{\mathcal{P}\mathrm{ol}_{#1}(\Sector,\Ef)}

\def\commutatif{\ar@{}[rd]|{\circlearrowleft}}

\newtheorem{theorem}{Theorem}[section]

\newtheorem{proposition}[theorem]{Proposition}
\newtheorem{lemma}[theorem]{Lemma}

\newtheorem{innercustomthm}{Theorem}
\newenvironment{customthm}[1]
{\renewcommand\theinnercustomthm{#1}\innercustomthm}
{\endinnercustomthm}

\theoremstyle{definition}

\newtheorem{definition}[theorem]{Definition}
\newtheorem{remark}[theorem]{Remark}

\theoremstyle{definition}

\theoremstyle{definition}
\newtheorem{notation}[theorem]{Notation}

\usepackage{bbold}

\usepackage{natbib}
\usepackage{graphicx}
\usepackage[all]{xy}
\usepackage{xypic}
\usepackage{amssymb}
\usepackage{amsmath}
\usepackage{mathrsfs}
\usepackage{mathtools}
\usepackage{amsthm}
\usepackage{moreverb}
\usepackage{geometry}

\title{POLYNOMIAL FUNCTORS ON SOME CATEGORIES OF ELEMENTS}
\author{OURIEL BLOEDE}

\date{}
\begin{document}

\maketitle

\begin{abstract}

We study the category $\Funcat$ of functors from the category $\Sector$, which is the category of elements of some presheaf $S$ on the category $\E$ of finite dimensional vector spaces, to $\Ef$ the category of vector spaces of any dimension on some field $\F$.

In the case where $S$ satisfies some noetherianity condition, we have a convenient description of the category $\Sector$. In this case, we can define a notion of polynomial functors on $\Sector$. And, like in the usual setting of functors from the category of finite dimensional vector spaces to the one of vector spaces of any dimension, we can describe the quotient $\pol{n}/\pol{n-1}$, where $\pol{n}$ denote the full subcategory of $\Funcat$ of polynomial functors of degree less than or equal to $n$.

Finally, if $\F=\Fp$ for some prime $p$ and if $S$ satisfies the required noetherianity condition, we can compute the set of isomorphism classes of simple objects in $\Funcat$.

\end{abstract}

\section{Introduction}

\begin{notation}
 In the following, for $F$ a functor on a category $\mathcal{C}$, $c$ and $c'$ objects of $\mathcal{C}$ and $f\ :\ c\rightarrow c'$ an arrow in $\mathcal{C}$, if there is no ambiguity on $F$, we will denote by $f_*$ the induced map $F(f)$ from $F(c)$ to $F(c')$ if $F$ is covariant, and by $f^*$ the induced map from $F(c')$ to $F(c)$ if $F$ is contravariant.\\ 
 
 We denote by $\SVF$  the category of contravariant functors from $\E$ the category of finite dimensional vector spaces over a given field $\F$ to $\mathcal{S}\mathrm{et}$ the category of sets, and by $\FVF$ the full subcategory of $\SVF$ with objects the functors with values in $\mathcal{F}\mathrm{in}$ the category of finite sets.
\end{notation}

\subsection{The category $\Funcat$}

For $S\in\SVF$, the category of elements $\Sector$ is the category whose objects are the pairs $(W,\psi)$ with $W\in\E$ and $\psi$ in $S(W)$ and whose morphisms from $(W,\psi)$ to $(H,\eta)$ are the morphisms $\gamma$ of $\F$-vector spaces from $W$ to $H$ satisfying $\gamma^*\eta=\psi$.\\ 

The aim of this article is to study the category $\Funcat$ of functors from the category $\Sector$ to the category $\Ef$ of $\F$-vector spaces of any dimension, under some conditions on $S$. \\

Our motivations to study such categories come from the theory of unstable modules over the Steenrod algebra. We explain succinctly (see \cite{B2}) how such categories appear in the study of unstable modules over an unstable algebra $K$ over the mod $2$ Steenrod algebra. For $K$ in $\K$ the category of unstable algebras, we consider the functor $S$ that maps the vector space $W$ to $\HomK(K,H^*(W))$, for $H^*(W)$ the cohomology with coefficients in $\mathbb{F}_2$ of the classifying space $BW$. This functor takes it values in the category of profinite sets. \\

The functor that maps $W$ to $\mathbb{F}_2^{S(W)}$ (the set of continuous maps from $S(W)$ to $\mathbb{F}_2$) is then an algebra in the category $\mathcal{F}(\E,\Ef)$ of functors from $\E$ to $\Ef$. For $K-\U$ the category of unstable $K$-modules, and $\Nil$ the localising subcategory of nilpotent modules, we have an equivalence of categories between $K-\U/\Nil$ and the full subcategory of analytic functors in $\mathbb{F}_2^S-\mathcal{M}od$.\\

In the case where $K$ is noetherian, $S$ takes values in $\mathcal{F}in$ and $\mathbb{F}_2^S-\mathcal{M}od\cong\Funcat$ (cf \cite{djament:tel-00119635}).\\

Since any simple object in $K-\U$ is the suspension of a $nil$-closed simple object in $K-\U$, the computation of simple objects in $\Funcat$ would allow us to classify simple objects in $K-\U$.\\

In section \ref{partintr}, we recall the definition of the kernel of an element of $\SVF$ as well as the definition of a noetherian functor from \cite{HLS2}. We use those to describe the category $\Sector$ in the case where $S$ satisfies a condition slightly weaker than the noetherianity condition of \cite{HLS2}.

\subsection{Polynomial functors in $\Funcat$}

In section \ref{part2}, we define and study a notion of polynomial functor in $\mathcal{F}(\Sector,\E)$. Polynomial functors over an additive category are already well studied and have very interesting properties such as homological finiteness. They have been of importance in computing the simple objects of the category $\mathcal{F}(\E,\Ef)$ (see for example \cite{PS}). The category $\Sector$ is not additive, yet in the case where $S$ satisfies the weaker noetherianity condition, we can still introduce a notion of polynomiality. For $\pol{n}$ the full subcategory of polynomial functors of degree $n$ on $\Sector$, we get the following theorem, where the category $\mathfrak{R}_S$, that we will introduce in the first section, is equivalent to a category with a finite set of objects, in the case where $S$ is noetherian in the sense of \cite{HLS2}.

\begin{customthm}{\ref{main1}}
    There is an equivalence of categories between $\pol{n}/\pol{n-1}$ and\\ $\mathcal{F}(\mathfrak{R}_S,\F\left[\Sigma_n\right]-\mathcal{M}\mathrm{od})$.
\end{customthm}

\subsection{Simple functors in $\Funcat$}

In the case where $\F$ is a finite field $\Fp$ with $p$ prime, using similar techniques to those presented in \cite{PS}, we are able to describe simple objects in $\Funcat$.

\begin{customthm}{\ref{mainx}}
    There is a one-to-one correspondence between isomorphism classes of simple objects of $\Funcat$ and isomorphism classes of simple objects of $$\bigsqcup\limits_{(W,\psi),n}\Fp\left[\mathcal{A}ut_{\Sector}(W,\psi)\times\Sigma_n\right]-\mathcal{M}od$$ with $(W,\psi)$ running through the isomorphism classes of objects in $\mathfrak{R}_S$ and $n$ running through $\N$.
\end{customthm}

\textbf{Acknowledgements: }I am thankful to Geoffrey Powell for his careful proofreading and for his continued support during and after my PhD. This work has been partially supported by the Labex CEMPI (ANR-11-LABX-0007-01).

\section{Noetherian functors}\label{partintr}

In this section, we start by recalling the definition of a noetherian functor from \cite{HLS2} and we introduce the weaker noetherianity condition that will be needed in the following sections.\\

For $S$ satisfying the weaker noetherianity condition, the category $\Sector$ can be described using Rector's category $\Rector$. We introduce this category, and end this section by comparing the categories of functor $\Funcat$ and $\mathcal{F}(\Rector\times\E,\Ef)$. 

\subsection{Definition and first properties}
We start by recalling the definition of the kernel of an element of $S(V)$, for $S\in\SVF$ and $V\in\E$.

\begin{proposition}\cite[Proposition-Definition 5.1]{HLS2}\label{imp}
Let $S\in\SVF$, $V\in\E$ and $s\in S(V)$. Then, there exists a unique sub-vector space $U$ of $V$, denoted by $\text{ker}(s)$, such that: 
\begin{enumerate}
    \item For all $t\in S(W)$ and all morphism $\alpha\ :\ V\rightarrow W$ such that $s=\alpha^*t$, $\text{ker}(\alpha)\subset U$.
    \item There exists $W_0$ in $\E$, $t_0\in S(W_0)$ and $\alpha_0\ :\ V\rightarrow W_0$ such that $s=\alpha_0^*t$ and $\text{ker}(\alpha_0)=U$.
    \item There exists $t_0\in S(V/U)$ such that $s=\pi^*t_0$, where $\pi$ is the projection of $V$ onto $V/U$.
\end{enumerate}
\end{proposition}

Notice that, since $\pi\ :\ W\rightarrow W/\ker(\psi)$ is surjective, it admits a right inverse, therefore $\pi^*$ has a left inverse, hence it is injective. We will denote by $\Tilde{\psi}$ the unique element of $S(W/\ker(\psi))$ such that $\pi^*\Tilde{\psi}=\psi$.

\begin{definition}
Let $S\in\SVF$, $V\in\E$ and $s\in S(V)$. We say that $s$ is regular if $\ker(s)=0$. 

Let $\text{reg}(S)(V):=\{x\in S(V)\ ;\ \ker(x)=0\}$.
\end{definition}

We also recall the definition of Rector's category $\mathfrak{R}_S$ which is the full subcategory of $\Sector$ whose objects are the pairs $(W,\psi)$ with $\psi$ regular.\\

\begin{definition}\label{noetherianfunc}
Let $S$ be in $\FVF$, we say that $S$ is noetherian if it satisfies the following:

\begin{enumerate}
    \item there exists an integer $d$ such that $\text{reg}(S)(V)=\emptyset$ for $\dim(V)> d$,
    \item\label{noethcond} for all $V\in\E$ and $s\in S(V)$ and for all morphisms $\alpha$ which takes values in $V$, $\text{ker}(\alpha^*s)=\alpha^{-1}(\text{ker}(s))$.
\end{enumerate} 
\end{definition}

In \cite{HLS2}, the authors proved that $S\in\FVF$ is noetherian if and only if there is a noetherian unstable algebra $K$ such that $S\cong\HomK(K,H^*(\_))$. In this article, $S$ will not need to satisfy all conditions of Definition \ref{noetherianfunc}.

\begin{definition}
We say that $S\in\SVF$ satisfies the weaker noetherianity condition if it satisfies condition \ref{noethcond} in Definition \ref{noetherianfunc}.
\end{definition} 

Yet, our results will be of particular interest in the case where $S$ is noetherian, since in this case Rector's category admits a finite skeleton.\\

\begin{definition}
    An object $S\in\SVF$ is connected if $S(0)$ has a single element $\epsilon$. In this case, for $V$ an object in $\E$, $\epsilon_V:=\pi_0^{V*}\epsilon$ for $\pi_0^V$ the unique map from $V$ to $0$.
\end{definition}

\begin{remark}
    In the case where $S$ is not connected, for $\gamma\in S(0)$, we can consider $S^\gamma$ that maps $W$ to the set of elements $\psi\in S(W)$ such that $0^*\psi=\gamma$. $S^\gamma$ is then a subfunctor of $S$ and $(S^\gamma)_{\gamma\in S(0)}$ is a partition of $S$. We get that $\Sector$ is the coproduct of the categories $\mathfrak{S}_{S^\gamma}$ and that a functor on $\Sector$ is just a family of functors over each of the categories $\mathfrak{S}_{S^\gamma}$.
\end{remark}

\subsection{The category $\Sector$}

In this subsection we describe the objects and morphisms in the category $\Sector$ in the case where $S$ is connected and satisfies the weaker noetherianity condition.

\begin{proposition}\label{defbox}
    We consider $S\in\SVF$ connected that satisfies the weaker noetherianity condition. Then, for any $(W,\psi)\in\Sector$, there exists a unique element $\psi\boxplus\epsilon_V\in S(W\oplus V)$, such that $\iota_W^*\psi\boxplus\epsilon_V=\psi$ and $\iota_V^*\psi\boxplus\epsilon_V=\epsilon_V$, for $\iota_W$ and $\iota_V$ the inclusions of $W$ and $V$ in $W\oplus V$.
\end{proposition}

\begin{proof}
    Let $\psi$ in $S(W)$. We consider $\pi^*\psi\in S(W\oplus V)$ for $\pi$ the projection from $W\oplus V$ to $W$ along $V$. It satisfies $\iota_W^*\pi^*\psi=\psi$. Furthermore, $\pi\circ\iota_V=0$. Hence, $\iota_V^*\pi^*\psi=0^*\psi$. Since the trivial morphism from $V$ to $W$ factorizes through the trivial vector space $0$, and since $S(0)=\{\epsilon\}$, $0^*\psi=0^*\epsilon$, $\iota_V^*\pi^*\psi=\epsilon_V$. This proves the existence condition. We now prove the uniqueness.\\

    For $\gamma\in S(W\oplus V)$ such that $\iota_W^*\gamma=\psi$ and $\iota_V^*\gamma=\epsilon_V$, since $S$ satisfies the weaker noetherianity condition, $V=\ker(\iota_V^*\gamma)=\iota_V^{-1}(\ker(\gamma))$. Therefore, $V\subset\ker(\gamma)$. By definition of the kernel, there exists $\Tilde{\gamma}\in S(W)$ such that $\gamma=\pi^*\Tilde{\gamma}$. Then, $\psi=\iota_W^*\pi^*\Tilde{\gamma}$, since $\pi\circ\iota_W=\id_W$, $\psi=\Tilde{\gamma}$. Which prove the uniqueness condition.
\end{proof}

The notation $\psi\boxplus\epsilon_V$ will be convenient in the following, but as we have seen, it is just the element $\pi^*\psi\in S(W\oplus V)$ for $\pi$ the projection from $W\oplus V$ onto $W$.\\

By definition of the kernel, for any $W\in\E$ and $\psi\in S(W)$, there exists a unique $\Tilde{\psi}\in S(W/\ker(\psi))$ such that $$\psi=\pi^*\Tilde{\psi}=\Tilde{\psi}\boxplus\epsilon_{\ker(\psi)}.$$ Since $S$ satisfies the weaker noetherianity condition, $\Tilde{\psi}$ is regular (this is because $\ker(\psi)=\pi^{-1}(\ker(\Tilde{\psi}))$). We get the following lemma.

\begin{lemma}\label{2.8.2}
    For $S$ connected that satisfies the weaker noetherianity condition and for $(W,\psi)\in\Sector$, $(W,\psi)\cong(W/\ker(\psi)\oplus\ker(\psi),\Tilde{\psi}\boxplus\epsilon_{\ker(\psi)})$, with $(W/\ker(\psi),\Tilde{\psi})\in\Rector$.
\end{lemma}

We now describe morphisms in $\Sector$, using this decomposition.

\begin{proposition}\label{2.8}
    Let $(W,\psi)$ and $(H,\eta)$ be two objects in $\mathfrak{R}_S$, and let $U$ and $V$ be two finite dimensional vector spaces. The set of morphisms in $\Sector$ from $(W\oplus U,\psi\boxplus\epsilon_U)$ to $(H\oplus V,\eta\boxplus\epsilon_V)$ is the set of morphisms $\alpha$ whose block matrices have the form $\left(\begin{array}{cc}
    f & 0 \\
    g & h
\end{array}\right)$, with $f$ a morphism from $(W,\psi)$ to $(H,\eta)$ in $\Rector$, $g$ a morphism from $W$ to $V$ and $h$ a morphism from $U$ to $V$.
\end{proposition} 

\begin{proof}
    First, we prove that for such $\alpha$, $\alpha$ is a morphism from $(W\oplus U,\psi\boxplus\epsilon_U)$ to $(H\oplus V,\eta\boxplus\epsilon_V)$ in $\Sector$. We have $\iota_W^*\alpha^*(\eta\boxplus\epsilon_V)=\iota_W^*\alpha^*\pi^*(\eta)$ for $\pi$ the projection from $H\oplus V$ onto $H$. This is equal to $f^*\eta=\psi$. Also, $\iota_U^*\alpha^*(\eta\boxplus\epsilon_V=h^*\epsilon_V=\epsilon_U$. Then, by Proposition \ref{defbox}, $\alpha^*(\eta\boxplus\epsilon_V)=\psi\boxplus\epsilon_U$.\\

    We now prove that morphisms from $(W\oplus U,\psi\boxplus\epsilon_U)$ to $(H\oplus V,\eta\boxplus\epsilon_V)$ have this form. First, we have that $$U=\ker(\psi\boxplus\epsilon_U)=\alpha^{-1}(\ker(\eta\boxplus\epsilon_V))=\alpha^{-1}(V).$$ Hence, $\alpha(U)\subset V$. Now, we consider the composition $\pi_H\circ\alpha$ from $(W\oplus U)$ to $H$ for $\pi_H$ the projection from $H\oplus V$ onto $H$. We have $\pi_H\circ\alpha=f\circ\pi_W$ for $\pi_W$ the projection from $(W\oplus U)$ onto $W$. Then, $\psi\boxplus\epsilon_U=\alpha^*\pi_H^*\eta=\pi_W^*(f^*\eta)$. We get, since $\pi_W^*$ is injective, that $f^*\eta=\psi$, therefore $f$ is a map from $(W,\psi)$ to $(H,\eta)$ in $\Rector$. This concludes the proof.
\end{proof}

\begin{remark}
    It is worth noticing that, since $S$ satisfies the weaker noetherianity condition, morphisms from $(W,\psi)$ to $(H,\eta)$ in $\Rector$ are necessarily injective morphisms from $W$ to $H$. This is one reason why functors on $\Rector$ are a lot easier to understand than functors on $\Sector$, and it will be a key fact in computing simple objects in $\Funcat$.
\end{remark}

\subsection{The categories $\Funcat$ and $\mathcal{F}(\Rector\times\E,\Ef)$}

In this subsection, we compare the categories $\Funcat$ and $\mathcal{F}(\Rector\times\E,\Ef)$, with $S$ connected and satisfying the weaker noetherianity condition.\\

By Lemma \ref{2.8.2}, for $W\in\E$ and $\psi\in S(W)$, $(W,\psi)$ is isomorphic as an object of $\Sector$ with $(W/\ker(\psi)\oplus\ker(\psi),\Tilde{\psi}\boxplus\epsilon_{\ker(\psi)})$. Therefore, we have a faithfull and essentially surjective functor from $\mathfrak{R}_S\times\E$ to $\Sector$ that maps the pair $((W,\psi),V)$ with $\psi$ regular to $(W\oplus V,\psi\boxplus\epsilon_V)$.\\

This functor is not full, indeed (Proposition \ref{2.8}) the set of morphisms between $(W\oplus V,\psi\boxplus\epsilon_V)$ and $(H\oplus U,\eta\boxplus\epsilon_U)$ in $\Sector$ is given by the linear maps whose block matrices are of the form $\left(\begin{array}{cc}
    f & 0 \\
    g & h
\end{array}\right)$ with $g$ and $h$ any linear maps respectively from $W$ and $\ker(\psi)$ onto $\ker(\eta)$ and $f$ a morphism in $\mathfrak{R}_S$ from $(W,\psi)$ to $(H,\eta)$, whereas the image of  $\mathfrak{R}_S\times\E$ contains only maps of the form $\left(\begin{array}{cc}
    f & 0 \\
    0 & h
\end{array}\right)$. Yet, it admits a left quasi-inverse that maps $(W,\psi)$ to $((W/\ker(\psi),\Tilde{\psi}),\ker(\psi))$ which is full and essentially surjective but not faithful. More precisely, two maps from $(W,\psi)$ to $(H,\eta)$ have the same image if and only if their restriction to $\ker(\psi)$ are equal as well as their induced maps from $W/\ker(\psi)$ to $H/\ker(\eta)$.

\begin{definition}\label{defi}
    Let $\mathcal{O}$ be the functor from $\Funcat$ to $\mathcal{F}(\Rector\times\E,\Ef)$ induced by the functor from $\Rector\times\E$ to $\Sector$ that maps $((W,\psi),V)\in\Rector\times\E$ to $(W\oplus V,\psi\boxplus\epsilon_V)$ and the morphism $(f,h)$ in $\Rector\times\E$ to $\left(\begin{array}{cc}
    f & 0 \\
    0 & h
\end{array}\right)$ in $\Sector$.\\

Let also $\mathcal{E}$ be the functor from $\mathcal{F}(\Rector\times\E,\Ef)$ to $\Funcat$ induced by the functor that maps $(W,\psi)$ in $\Sector$ to $((W/\ker(\psi),\Tilde{\psi}),\ker(\psi))$ in $\Rector\times\E$ and $f$ from $(W,\psi)$ to $(H,\eta)$ to $(\Tilde{f},f|_{\ker(\psi)})$ with $\Tilde{f}$ the morphism induced by $f$ from $(W/\ker(\psi),\Tilde{\psi})$ to $(H/\ker(\eta),\Tilde{\eta})$.
\end{definition}

\begin{lemma}\label{2zoug}
    For $\lambda$ a natural transformation in $\mathcal{F}(\Rector\times\E,\Ef)$ from $G$ to $\mathcal{O}(F)$, $\lambda$ extends to a natural transformation in $\Funcat$ from $\mathcal{E}(G)$ to $F$ if and only if, for any $(W,\psi)\in\Rector$, $V\in\E$ and $f\in\HomF(W,V)$, the following diagram commutes:
    $$\xymatrix{ & F(W\oplus V,\psi\boxplus\epsilon_V)\ar[dd]^-{\alpha_*}\\
    G((W,\psi),V)\ar[ru]^{\lambda}\ar[rd]_-\lambda &\\
    & F(W\oplus V,\psi\boxplus\epsilon_V),}$$
    with $\alpha$ the morphism whose block matrix is given by $\left(\begin{array}{cc}
    \id_W & 0 \\
    f & \id_V
\end{array}\right)$. 
\end{lemma}

\begin{proof}
    The only if part is straightforward. Let's assume that $\lambda$ satisfies the required condition. Then, for $(W,\psi)$ in $\Sector$, one can choose arbitrarily a complementary subspace $C$ of $\ker(\psi)$, then for $\gamma$ the inverse isomorphism of the projection from $C$ to $W/\ker(\psi)$, one can define $\lambda_{(W,\psi)}$ from $G((W/\ker(\psi),\Tilde{\psi}),\ker(\psi))$ to $F(W,\psi)$ as the composition of $\lambda_{((W/\ker(\psi),\Tilde{\psi}),\ker(\psi))}$, which values are in $F(W/\ker(\psi)\oplus\ker(\psi),\Tilde{\psi}\boxplus\epsilon_{\ker(\psi)})$, with $$(\gamma\oplus\id_{\ker(\psi)})_*\ :\ F(W/\ker(\psi)\oplus\ker(\psi),\Tilde{\psi}\boxplus\epsilon_{\ker(\psi)})\rightarrow F(W,\psi).$$ 

    The required condition guarantees that this does not depend on the choice of $C$. Furthermore, it entails that $\lambda$ is a natural transformation on $\Sector$, since any morphism in $\Sector$ from $(H,\eta)$ to $(W,\psi)$ can be factorised as $\left(\begin{array}{cc}
    \id_C & 0 \\
    f & \id_{\ker(\psi)}
\end{array}\right)\circ \left(\begin{array}{cc}
    g & 0 \\
    0 & h
\end{array}\right)$ with some morphisms $f$ and $h$ and some injective morphism $g$.
\end{proof}

\section{Polynomial functors over $\Sector$}\label{part2}

Since $\mathcal{F}(\mathfrak{R}_S\times\E,\Ef)$ is isomorphic to  $\mathcal{F}(\mathfrak{R}_S,\mathcal{F}(\E,\Ef))$, there is a notion of polynomial functors of degree $n$ for functors in $\mathcal{F}(\mathfrak{R}_S\times\E,\Ef)$ corresponding to those taking values in polynomial functors of degree $n$ from $\E$ to $\Ef$, in the sense of \cite{P1} or \cite{DV}. We denote by $\mathcal{P}\mathrm{ol}_n(\Rector\times\E,\Ef)$ the full subcategory of $\mathcal{F}(\mathfrak{R}_S\times\E,\Ef)$ of polynomial functors of degree less than or equal to $n$. Using purely formal arguments, as well as known facts about polynomial functors in $\mathcal{F}(\E,\Ef)$, one could easily compute the categorical quotient $\mathcal{P}\mathrm{ol}_n(\Rector\times\E,\Ef)/\mathcal{P}\mathrm{ol}_{n-1}(\Rector\times\E,\Ef)$ and would find that it is equivalent to $\mathcal{F}(\Rector,\F\left[\Sigma_n\right]-\mathcal{M}\mathrm{od})$. All the difficulties in the following section come from the fact that the functor from $\mathfrak{R}_S\times\E$ to $\Sector$ is not full. \\

In this section we define a notion of polynomial functors on $\Sector$ and manage to compute the quotient $\pol{n}/\pol{n-1}$.

\subsection{Definition}

We recall that for $F\in\mathcal{F}(\E,\Ef)$, $\Bar{\Delta} F(W)$ is the kernel of the map from $F(W\oplus\F)$ to $F(W)$ induced by the projection along $\F$. Polynomial functors of degree at most $n$ are functors $F$ such that $\Bar{\Delta}^{n+1}F=0$ and $\mathcal{P}\mathrm{ol}_n(\E,\Ef)$ denote the full subcategory of polynomial functors of degree at most $n$ in $\mathcal{F}(\E,\Ef)$. We define similar notions for functors on $\Sector$.\\

We start this section by defining a difference functor $\decal$ on $\Funcat$.\\

\begin{definition}
    $\Bar{\Delta}_{(\F,\epsilon_{\F})}\ :\ \mathcal{F}(\mathfrak{S}_S,\Ef)\rightarrow \mathcal{F}(\mathfrak{S}_S,\Ef)$ is the functor such that $\decal F(W,\psi)$ is the kernel of the map $F(W\oplus\F,\psi\boxplus\epsilon_{\F})\rightarrow  F(W,\psi)$ induced by the projection from $W\oplus\F$ to $W$, and such that, for $\alpha $ a morphism in $\Sector$, $\decal F(\alpha)$ is the map induced by $\alpha\oplus\id_{\F}$.
\end{definition}

\begin{lemma}\label{decexact}
    The functor $\decal$ is exact.
\end{lemma}

\begin{proof}
    We consider the following short exact sequence in $\Funcat$ :
    $$0\rightarrow F'\rightarrow F\rightarrow F"\rightarrow 0.$$
For $(W,\psi)\in\Sector$, we get the following commutative diagram whose lines are exact : 
$$\xymatrix{0\ar[r] & F'(W\oplus\F,\psi\boxplus\epsilon_{\F})\ar[r]\ar[d] & F(W\oplus\F,\psi\boxplus\epsilon_{\F})\ar[r]\ar[d] & F"(W\oplus\F,\psi\boxplus\epsilon_{\F})\ar[r]\ar[d] & 0\\
0\ar[r] & F'(W,\psi)\ar[r] & F(W,\psi)\ar[r] & F"(W,\psi)\ar[r] & 0,
}$$
whose vertical maps are induced by the projection from $W\oplus\F$ to $W$. Using the exactness of the lines and commutativity of the diagram, one checks that it induces a short exact sequence $$\xymatrix{
0\ar[r] & \decal F'(W,\psi)\ar[r] & \decal F(W,\psi)\ar[r] & \decal F"(W,\psi)\ar[r] & 0.
}$$
This exact sequence is natural in $(W,\psi)$.
\end{proof}

\begin{definition}\label{polanal}
    $F\in\Funcat$ is polynomial of degree less than $n$ if $\decal^{n+1}F=0$. We denote by $\pol{n}$ the full subcategory of $\Funcat$ whose objects are the polynomial functors of degree less than or equal to $n$.
\end{definition}

\begin{proposition}
    The category $\pol{n}$ is a Serre class of $\pol{n+1}$.
\end{proposition}

\begin{proof}
    This is straightforward from Lemma \ref{decexact}.
\end{proof}

There is a notion of analytic funtors on $\mathcal{F}(\E,\Ef)$. Those are the functors which are the colimit of their polynomial sub-functors. Similarly, one can define a notion of analytic functors on $\Sector$.\\

\begin{lemma}
    Let $F\in\Funcat$. $F$ admits a greatest polynomial sub-functor of degree less than or equal to $n$. We denote it by $p_n(F)$.
\end{lemma} 

\begin{proof}
    For $(W,\psi)\in\Sector$ and $x\in F(W,\psi)$, we denote by $<x>_F$ the image of $\F\left[\Hom_{\Sector}((W,\psi),(\_,\_))\right]$ under the natural morphism  that maps $\id_{(W,\psi)}$ to $x$. We say that $x$ is polynomial of degree less than or equal to $n$ if and only if $<x>_F$ is.\\

    The functor $F$ is polynomial of degree less than or equal to $n$ if and only if $x$ is polynomial of degree less than or equal to $n$ for any $(W,\psi)\in\Sector$ and any $x\in F(W,\psi)$. The condition is obviously necessary since $\decal^{n+1}<x>_F$ is a sub-functor of $\decal^{n+1}F$. If we assume that $F$ is not polynomial of degree less than or equal to $n$, $\decal^{n+1}F$ is not trivial, therefore there exists $(W,\psi)\in\Sector$ and $x\in\decal^{n+1}F(W,\psi)\subset F(W\oplus\F^{n+1},\psi\boxplus\epsilon_{\F^{n+1}})$ different from $0$. Then, $x$ is not polynomial of degree less than or equal to $n$. The condition is therefore sufficient.\\

    Finally, the set of elements $x$ of $F$ polynomial of degree less than or equal to $n$ defines a sub-functor of $F$ and is greater than any polynomial sub-functor of $F$ of degree less than or equal to $n$.
\end{proof}

By definition, we have $p_n(F)\subset p_{n+1}F$.

\begin{definition}
    A functor $F$ in $\Funcat$ is said to be analytic if it is the colimit of the $p_n(F)$.\\

    $\mathcal{F}_\omega(\Sector,\Ef)$ is the full subcategory of $\Funcat$ of analytic functors.
\end{definition}

   We end this subsection by computing $\pol{0}$, the following subsections will have the purpose of describing the quotients $\pol{n}/\pol{n-1}$

\begin{proposition}\label{pol0}
    The categories $\pol{0}$ and $\mathcal{F}(\mathfrak{R}_S,\Ef)$ are equivalent.
\end{proposition}

\begin{proof}
    For $F\in\pol{0}$, $(W,\psi)$ an object of $\Sector$ and $\Tilde{\psi}\in F(W/\ker(\psi))$ such that $\pi^*\Tilde{\psi}=\psi$, $\pi_*$ (induced by $\pi$ the projection in $\Sector$ from $(W,\psi)$ to $(W/\ker(\psi),\Tilde{\psi})$) is a natural isomorphism between $F(W,\pi^*\psi)$ and $F(W/\ker(\psi),\psi)$. Indeed, $\pi_*$ may be factorised in the following way $$F(W,\psi)\cong F(W/\ker(\psi)\oplus\F^k,\Tilde{\psi}\boxplus\epsilon_{\F^k})\rightarrow ...\rightarrow F(W/\ker(\psi)\oplus\F,\Tilde{\psi}\boxplus\epsilon_{\F})\rightarrow F(W/\ker(\psi),\Tilde{\psi}),$$ where $k$ is the dimension of $\ker(\psi)$ and each map is induced by the projection that omits the last factor $\F$. And since $\decal F=0$, each of those maps is an isomorphism.\\
    
    The forgetful functor from $\Funcat$ to $\mathcal{F}(\mathfrak{R}_S,\Ef)$ has a right quasi-inverse that maps $F\in\mathcal{F}(\mathfrak{R}_S,\Ef)$ to $\Bar{F}$, where $\Bar{F}(W,\psi):=F(W/\ker(\psi),\Tilde{\psi})$ and $\Bar{F}(\gamma)$, for $\gamma\ :\ (W,\psi)\rightarrow (H,\eta)$, is $F(\Tilde{\gamma})$ for $\Tilde{\gamma}$ the induced map from $(W/\ker(\psi),\Tilde{\psi})$ to $(H/\ker(\eta),\Tilde{\eta})$. By construction, it is a quasi-inverse, if we restraint the forgetful functor to $\pol{0}$. 
\end{proof}

\subsection{The $n$-th cross effect}

    In the context where $F$ is a functor over $\E$, the $n$-th cross effect $\mathrm{cr}_nF(X_1,...,X_n)$ is defined as the kernel of the map from $F(X_1\oplus...\oplus X_n)$ to $\bigoplus\limits_{1\leq i\leq n} F(X_1\oplus...\oplus \widehat{X_i}\oplus...\oplus X_n)$ induced by the projections from $X_1\oplus...\oplus X_n$ to $X_1\oplus...\oplus \widehat{X_i}\oplus...\oplus X_n$. Since $\Sigma_n$ acts on $\mathrm{cr}_nF(\F,...,\F)$ by permuting the factors $\F$, $F\mapsto\mathrm{cr}_nF(\F,...,\F)$ defines a functor from $\mathcal{F}(\E,\Ef)$ to $\F\left[\Sigma_n\right]-\mathcal{M}\mathrm{od}$.\\

    We consider $T^n$ the functor from $\E$ to itself that maps $V$ to $V^{\otimes n}$. $\Sigma_n$ has a right-action on $V^{\otimes n}$ with $v_1\otimes...\otimes v_n\cdot\sigma=v_{\sigma^{-1}(1)}\otimes...\otimes v_{\sigma^{-1}(n)}$. We get the following Proposition from \cite{P1}.

    \begin{proposition}\label{pirash}
    The functor from $\mathcal{P}\mathrm{ol}_n(\E,\Ef)$ to $\F\left[\Sigma_n\right]-\mathcal{M}\mathrm{od}$ that maps $F$ to $\mathrm{cr}_nF(\F,...,\F)$ is right adjoint to the functor that maps $M\in\F\left[\Sigma_n\right]-\mathcal{M}\mathrm{od}$ to $T^n\otimes_{\Sigma_n}M$.
        
    \end{proposition}

    Throughout this subsection $n$ is a fixed positive integer. We want to describe the quotient category $\pol{n}/\pol{n-1}$. To do so, we introduce a cross effect functor for functors on $\Sector$.\\

\begin{lemma}\label{kernuck}
    Let $F\in\Funcat$ and $(W,\psi)$ an object in $\Sector$. $\decal^n F(W,\psi)$ is the kernel of the map from $F(W\oplus\F^n,\psi\boxplus\epsilon_{\F^n})$ to $\bigoplus\limits_{i=1}^n F(W\oplus \F^{i-1}\oplus\widehat{\F}\oplus \F^{n-i},\psi\boxplus \epsilon_{\F^{n-1}})$ induced by the projections from $(W\oplus\F^n,\psi\boxplus\epsilon_{\F^n})$ to $(W\oplus \F^{i-1}\oplus\widehat{\F}\oplus \F^{n-i},\psi\boxplus \epsilon_{\F^{n-1}})$ in $\Sector$.
\end{lemma}

\begin{proof}
    This is straightforward by induction.
\end{proof}


More generally, we consider the n-th cross effect $\mathrm{cr}_nF$ defined as follows.

\begin{definition}
    For $F\in\Funcat$, $\mathrm{cr}_nF$ is the functor from $\Sector\times (\E)^n$ to $\Ef$ where $\mathrm{cr}_n(W,\psi;X_1,...,X_n)$ is the kernel of $F(W\oplus X_1\oplus...\oplus X_n,\psi\boxplus\epsilon)$ to $\bigoplus\limits_{i=1}^n F(W\oplus X_1\oplus...\oplus\widehat{X_i}\oplus...\oplus X_n,\psi\boxplus \epsilon)$ induced by the projections from $(W\oplus X_1\oplus...\oplus X_n,\psi\boxplus\epsilon)$ to $(W\oplus X_1\oplus...\oplus\widehat{X_i}\oplus ...\oplus X_n,\psi\boxplus \epsilon)$ in $\Sector$.
\end{definition}

We can restate Lemma \ref{kernuck} as $$\decal^n F(W,\psi)=\mathrm{cr}_nF(W,\psi;\F,...,\F).$$

As in the classical case of functors on $\E$, the $n$-th symmetric group acts on $\mathrm{cr}_nF(\_,\_;\F,...,\F)$ by permuting the factors $\F$. Therefore, $\mathrm{cr}_nF(\_,\_;\F,...,\F)$ takes values in $\F\left[\Sigma_n\right]-\mathcal{M}od$, the category of $\Sigma_n$-representations over $\F$.\\

Considering functors on $\Rector\times\E$, one can check easily that Proposition \ref{pirash} gives rise to a similar adjunction from $\mathcal{P}\mathrm{ol}_n(\Rector\times\E,\Ef)$ to $\mathcal{F}(\Rector,\F\left[\Sigma_n\right]-\mathcal{M}\mathrm{od})$ whose left adjoint maps $M$ to $T^n\otimes_{\Sigma_n}M$, that maps $((W,\psi),V)\in\Rector\times\E$ to $V^{\otimes n}\otimes_{\Sigma_n}M(W,\psi)$.\\

We want to extend this adjunction to $\Funcat$. To do so, we need to prove that $\mathrm{cr}_nF(\_,\_;\F,...,\F)$ behaves well with respect to the maps in $\Sector$ that are not obtained from maps in $\Rector\times\E$.\\

The end of this subsection is devoted to prove that if $F\in\pol{n}$ then $\mathrm{cr}_nF(\_,\_;\F,...,\F)$ does behave well with respect to those maps.\\

We notice that for $F$ polynomial of degree $n$, $\mathrm{cr}_nF(W,\psi;\_,...,\_)$ is additive in each variable. More explicitly:\\ 

\begin{lemma}\label{constanti}
    $\mathrm{cr}_nF(W,\psi;X_1\oplus Y_1,...,X_n\oplus Y_n)$ is isomorphic to $\bigoplus \mathrm{cr}_nF(W,\psi;A_1,...,A_n)$, with the direct sum going through the families $(A_1,...,A_n)$ with $A_i=X_i\text{ or }Y_i$. The isomorphism is given by the direct sum of the $\mathcal{Q}_{A_1,...,A_n}$ induced by the projections from $X_i\oplus Y_i$ onto $A_i$.
\end{lemma}

\begin{proof}
    The map $\mathcal{Q}_{A_1,...,A_n}$ is induced by the map from $\mathrm{cr}_nF(W,\psi;X_1\oplus Y_1,...,X_n\oplus Y_n)$ to $F(W\oplus A_1\oplus ...A_n,\psi\oplus \epsilon)$, with $A_i=X_i\text{ or }Y_i$, induced by the projection from $W\oplus (X_1\oplus Y_1)\oplus ...\oplus(X_n\oplus Y_n)$ to $W\oplus A_1\oplus ...\oplus A_n$.\\
    
    Since $\decal^{n+1}F=0$, it's kernel is precisely the direct sum of the images of the $\mathrm{cr}_nF(W,\psi;B_1,...,B_n)$ with at least one $B_i\neq A_i$, under the injections in $W\oplus (X_1\oplus Y_1)\oplus ...\oplus(X_n\oplus Y_n)$. The restriction to $\mathrm{cr}_nF(W,\psi;A_1,...,A_n)$ (seen as a subspace of $\mathrm{cr}_nF(W,\psi;X_1\oplus Y_1,...,X_n\oplus Y_n)$) is the identity.
\end{proof}

\begin{remark}
    The inverse of $\bigoplus\limits_{(A_1,...,A_n)}\mathcal{Q}_{A_1,...,A_n}$ from $\mathrm{cr}_nF(W,\psi;X_1\oplus Y_1,...,X_n\oplus Y_n)$ to\\ $\bigoplus \mathrm{cr}_nF(W,\psi;A_1,...,A_n)$ is given by the direct sum of the $\mathcal{I}_{A_1,...,A_n}$ which are the maps induced by the inclusions of the $A_i$ in $X_i\oplus Y_i$.\\

    We want to emphasize that the image of $\bigoplus \mathrm{cr}_nF(W,\psi;A_1,...,A_n)$ in $\bigoplus \mathrm{cr}_nF(W,\psi;U_1,...,U_n)$, for $A_i$ a sub-vector space of $U_i$ for each $i$, does not depend on the choice of complementary subspaces $B_i$, therefore the component of an element of $\mathrm{cr}_nF(W,\psi;U_1,...,U_n)$ in $\mathrm{cr}_nF(W,\psi;A_1,...,A_n)$ under the isomorphism of Lemma \ref{constanti} is the same for each choice of decomposition of the $U_i$ as $U_i=A_i\oplus B_i$. It will have some importance in the proof of Lemma \ref{ident}
\end{remark}

We consider $F\in\pol{n}$, $(W\oplus X_1\oplus...\oplus X_n,\psi\boxplus\epsilon)$ in $\Sector$ and a map $\alpha$ from $(W\oplus X_1\oplus...\oplus X_n,\psi\boxplus\epsilon)$ to itself of the form $\alpha=\left(\begin{array}{cc}
        \id_W & 0 \\
       f  & \id_{X_1\oplus...\oplus X_n}
    \end{array}\right)$. We have the following Lemma, with $\alpha_*$ the induced map from $F(W\oplus X_1\oplus...\oplus X_n,\psi\boxplus\epsilon)$.

\begin{lemma}\label{ident} 
$\alpha_*$ acts on $\mathrm{cr}_nF(W,\psi;X_1,...,X_n)$ as the identity.
\end{lemma}

\begin{proof}
  For $\pi_i$ that maps $x_1+...+x_n$ with $x_i\in X_i$ to $x_1+...+\hat{x}_i+...x_n$, we have that $\left(\begin{array}{cc}
        \id_W & 0 \\
       0  & \pi_i
    \end{array}\right)\circ\alpha$  is equal to  $$\left(\begin{array}{cc}
        \id_W & 0 \\
       \pi_i\circ f  & \pi_i
    \end{array}\right)= \left(\begin{array}{cc}
        \id_W & 0 \\
       \pi_i\circ f  & \id_{X_1\oplus...\oplus\widehat{X_i}\oplus...\oplus X_n}
    \end{array}\right)\circ\left(\begin{array}{cc}
        \id_W & 0 \\
       0  & \pi_i
    \end{array}\right).$$ Since $\mathrm{cr}_nF(W,\psi;X_1,...,X_n)$ is the intersection of the kernels of the $\left(\begin{array}{cc}
        \id_W & 0 \\
       0  & \pi_i
    \end{array}\right)_*$, this implies that the restriction of $\alpha_*$ to $\mathrm{cr}_nF(W,\psi;X_1,...,X_n)$ takes it values in $\mathrm{cr}_nF(W,\psi;X_1,...,X_n)$.\\

    We prove now that it acts as the identity. We consider the diagonal map $\Delta$ from $X_1\oplus...\oplus X_n$ to $(X_1\oplus X_1)\oplus...\oplus(X_n\oplus X_n)$ and the map $\alpha'$ from $W\oplus(X_1\oplus X_1)\oplus...\oplus(X_n\oplus X_n)$ to itself whose block matrix is given by $\left(\begin{array}{cc}
        \id_W & 0 \\
       \Delta\circ f  & \id_{X_1^{\oplus 2}\oplus...\oplus X_n^{\oplus 2}}
    \end{array}\right)$. It fits in the following commutative diagram:

\begin{center}
    \begin{tikzpicture}[scale=0.8, every node/.style={transform shape}]
\node[] (M1) at (7,7) {$F(W\oplus X_1\oplus...\oplus X_n,\psi\boxplus\epsilon)$};
\node[] (M2) at (7,4) {$F(W\oplus(X_1\oplus X_1)\oplus...(X_n\oplus X_n),\psi\boxplus\epsilon)$};
\node[] (M3) at (7,1) {$F(W\oplus(X_1\oplus X_1)\oplus...(X_n\oplus X_n),\psi\boxplus\epsilon)$};
\node[] (L3) at (1,1) {$F(W\oplus X_1\oplus...\oplus X_n,\psi\boxplus\epsilon)$};
\node[] (R3) at (13,1) {$F(W\oplus X_1\oplus...\oplus X_n,\psi\boxplus\epsilon),$};
\draw[->,>=latex] (M1) edge[in=90,out=-90] (M2);
\draw[->,>=latex] (M2) edge[in=90,out=-90] node[right]{$\alpha' _*$}(M3);
\draw[->,>=latex] (M1) edge[in=90,out=-156] node[left]{$\id_*$} (L3);
\draw[->,>=latex] (M1) edge[in=90,out=-32] node[right]{$\alpha_*$} (R3);
\draw[->,>=latex] (M3) edge[in=0,out=180] (L3);
\draw[->,>=latex] (M3) edge[in=180,out=0]  (R3);
        \end{tikzpicture}
\end{center}
where the top vertical map is induced by the injection of the $X_i$ in the first factor of $X_i\oplus X_i$, the right horizontal one is given by the projection of $X_i\oplus X_i$ on the first factor and the left horizontal one is the projection onto the first factor along the diagonal $\Delta(X_i)$ (i.e. the morphisms that maps $(x,y)$ in $X_i\oplus X_i$ to $x-y$ in $X_i$).\\ 

As we have seen, $\alpha'_*$ maps $\mathrm{cr}_nF(W,\psi,X_1\oplus X_1,...,X_n\oplus X_n)$ to itself. Also, the first factor $X_i$ of $X_i\oplus X_i$ admits two relevant complementary subspaces in $X_i\oplus X_i$. The first one is the second factor $X_i$, the second one is the diagonal $\Delta(X_i)$, since $X_i\oplus X_i=X_i\oplus\Delta(X_i)$ using that $(x,y)=(x-y,0)+(y,y)$ for $x$ and $y$ in $X_i$.\\ 

By Lemma \ref{constanti}, we get $$\mathrm{cr}_nF(W,\psi,X_1\oplus X_1,...,X_n\oplus X_n)\cong\bigoplus \mathrm{cr}_nF(W,\psi;A_1,...,A_n)\cong\bigoplus \mathrm{cr}_nF(W,\psi;A'_1,...,A'_n),$$ where the $A_i$ are either the first or the second factor in $X_i\oplus X_i$, and the $A'_i$ are either the first factor or the diagonal of $X_i\oplus X_i$.\\ 

The components $\mathrm{cr}_nF(W,\psi,X_1,...,X_n)$ where all $A_i$ and $A'_i$ are taken to be the first factor identifies under those isomorphism, and they are stable under $\alpha'_*$. From the left part of the commutative diagram above, we get that the restriction of $\alpha'_*$ to that component is the identity, which implies that $\alpha_*$ restricted to $\mathrm{cr}_nF(W,\psi;X_1,...,X_n)$ is also the identity. 
\end{proof}

\subsection{The category $\pol{n}/\pol{n-1}$}
In this subsection, we finally prove that $\decal^{n}$ induces an equivalence of categories between the localisation $\pol{n}/\pol{n-1}$ and the category $\mathcal{F}(\mathfrak{R}_S,\F\left[\Sigma_n\right]-\mathcal{M}\mathrm{od})$.\\ 

By abuse of notation, for $M$ a functor from $\Rector$ to $\F\left[\Sigma_n\right]-\mathcal{M}\mathrm{od}$, we denote by $T^n\otimes_{\Sigma_n}M$ the functor $\mathcal{E}(T^n\otimes_{\Sigma_n}M)$ (Definition \ref{defi}), which is the functor on $\Sector$ that maps $(W,\psi)$ to $\ker(\psi)^n\otimes_{\Sigma_n}M(W/\ker(\psi),\Tilde{\psi})$.

The following lemma is straightforward.

\begin{lemma}
     $T^n\otimes_{\Sigma_n}M\in \pol{n}$.
\end{lemma}

\begin{lemma}\label{congu}
    $\decal^n (T^n\otimes_{\Sigma_n}M)\cong M$ as a functor from $\mathfrak{R}_S$ to $\F\left[\Sigma_n\right]-\mathcal{M}\mathrm{od}$.
\end{lemma}

\begin{proof}
    From Lemma \ref{kernuck} and for $(W,\psi)$ regular, an element in $\decal^n(T^n\otimes_{\Sigma_n}M)(W,\psi)\subset T^n(\F^n)\otimes_{\Sigma_n} M(W,\psi)$ is mapped to $0$ under each map from $\F^n$ to $\F^{n-1}$ that send one of the factor $\F$ to $0$. Hence, an element of $\decal^n(T^n\otimes_{\Sigma_n}M)(W,\psi)$ admits a unique representing element in $T^n(\F^n)\otimes_{\F} M(W,\psi)$ of the form $v_1\otimes...\otimes v_n\otimes m$, with $(v_1,...,v_n)$ the canonical basis of $\F^n$ and $m\in M(W,\psi)$. We get the required isomorphism.
\end{proof}

\begin{proposition}
    The functor $M\mapsto T^n\otimes_{\Sigma_n}M$ is left adjoint to $$\decal^n\ :\ \pol{n}\rightarrow \mathcal{F}(\mathfrak{R}_S,\F\left[\Sigma_n\right]-\mathcal{M}\mathrm{od}).$$
\end{proposition}

\begin{proof}
    By naturality, a natural transformation from $T^n\otimes_{\Sigma_n}M$ to $F\in\Funcat$ is fully determined by the image of the class (for the equivalence relation induced by the actions of $\Sigma_n$ on $T^n$ and $M$) of the elements of the form $v_1\otimes...\otimes v_n\otimes m$ with $(W,\psi)$ an object of $\mathfrak{R}_S$, $m\in M(W,\psi)$ and  $(v_1,...,v_n)$ the canonical basis of $\F^n$. Furthermore, since $v_1\otimes...\otimes v_n\otimes m$ represents an element in $\decal^n(T^n\otimes_{\Sigma_n} M)(W,\psi)$, its image must be in $\decal^n F(W,\psi)$. Hence, the application that maps a natural transformation from $T^n\otimes_{\Sigma_n}M$ to $F$ to the induced morphism from $M$ to $\decal^{n}F$ is an injection.\\ 
    
    We consider a morphism in $\mathcal{F}(\mathfrak{R}_S,\F\left[\Sigma_n\right]-\mathcal{M}\mathrm{od})$ from $M$ to $\decal^nF$. By Proposition \ref{pirash}, it induces by adjunction a natural transformation of functors on $\Rector\times\E$ from $T^n\otimes_{\Sigma_n} M$ to $\mathcal{O}(F)$. Finally, Lemma \ref{ident} and Lemma \ref{2zoug} imply that, if $F$ is polynomial of degree $n$, this natural transformation can be extended as a natural transformation from $\mathcal{E}(T^n\otimes_{\Sigma_n} M)$ to $F$ in $\Funcat$.
\end{proof}

\begin{theorem}\label{main1}
    $\decal^n$ induces an equivalence of categories between $\pol{n}/\pol{n-1}$ and $\mathcal{F}(\mathfrak{R}_S,\F\left[\Sigma_n\right]-\mathcal{M}\mathrm{od})$.
\end{theorem}

\begin{proof}
    $\decal^n$ is an exact functor from $\pol{n}$ to $\mathcal{F}(\mathfrak{R}_S,\F\left[\Sigma_n\right]-\mathcal{M}\mathrm{od})$. It maps $\pol{n-1}$ to $0$, hence it induces a functor from $\pol{n}/\pol{n-1}$ to $\mathcal{F}(\mathfrak{R}_S,\F\left[\Sigma_n\right]-\mathcal{M}\mathrm{od})$. We consider $F\in\pol{n}$ and the exact sequence $0\rightarrow\ker(\eta)\rightarrow T^n\otimes_{\Sigma_n}\decal^n F\overset{\eta}{\rightarrow} F\rightarrow \coker(\eta)\rightarrow 0$, for $\eta$ the counit of the adjunction. By Lemma \ref{congu}, when we apply to it the functor $\decal^n$, the middle map becomes an isomorphism. Therefore, $\ker(\eta)$ and $\coker(\eta)$ are in $\pol{n-1}$, so $\eta$ is an isomorphism in $\pol{n}/\pol{n-1}$. We get that the functor induced by $\decal^n$ from $\pol{n}/\pol{n-1}$ to $\mathcal{F}(\mathfrak{R}_S,\F\left[\Sigma_n\right]-\mathcal{M}\mathrm{od})$ and the composition of $M\mapsto T^n\otimes_{\Sigma_n}M$ with the localization functor from $\pol{n}$ to $\pol{n}/\pol{n-1}$ are inverses.
\end{proof}

\section{Simple objects in $\Funcat$}

In this section, we describe the simple objects of the category $\Funcat$ for $\F$ a finite field $\Fp$ with $p$ prime, using the equivalence between $\pol{n}/\pol{n-1}$ and the category $\mathcal{F}(\mathfrak{R}_S,\Fp\left[\Sigma_n\right]-\mathcal{M}\mathrm{od})$. First, we prove that simple objects of $\Funcat$ are polynomial.\\

We consider the family of injective cogenerators $I_{(W\oplus V,\psi\boxplus \epsilon_V)}:=\Fp^{\Hom_{\Sector}(\_, (W\oplus V,\psi\boxplus\epsilon_V))}$.

\begin{proposition}
    For any $(W,\psi)\in \mathfrak{R}_S$ and any $V\in\E$, $\IWV$ is analytic.
\end{proposition}

\begin{proof}
    We have a forgetful functor from $\Sector$ to $\E$, it induces a functor $\mathcal{U}$ from $\mathcal{F}(\E,\Ef)$ to $\Funcat$. For $\Bar{\Delta}$ the usual difference functor in $\mathcal{F}(\E,\Ef)$, $\decal \mathcal{U}(F)(H,\eta)\cong \Bar{\Delta} F(H)$. Therefore,  $\decal^{n+1}\mathcal{U}(F)=0$ if and only if $\Bar{\Delta}^{n+1}F=0$, so $\mathcal{U}(F)\in\pol{n}$ if and only if $F$ is polynomial in the usual sense, and $F$ analytic implies $\mathcal{U}(F)$ analytic.\\

    By Proposition \ref{2.8}, $\Hom_{\Sector}((H,\eta),(W\oplus V,\psi\boxplus\epsilon_V))\cong\Hom_{\Sector}((H,\eta),(W,\psi)) \times\HomF(H,V)$. Therefore, $\IWV(H,\eta)$ is naturally isomorphic to the tensor product $$\Fp^{\Hom_{\Sector}((H,\eta),(W,\psi))}\otimes\Fp^{\HomF(H,V)}.$$
    We get that $\IWV\cong I_{(W,\psi)}\otimes \mathcal{U}(I_V)$, where $I_V$ denote the injective object in $\mathcal{F}(\E,\Ef)$ that maps $H$ to $\Fp^{\HomF(H,V)}$.\\ 
    
    $I_{(W,\psi)}$ is polynomial of degree $0$, indeed, since $(W,\psi)$ is regular we have that for any map from $(H\oplus\Fp,\eta\oplus\epsilon_{\Fp})$ to $(W,\psi)$ in $\Sector,$ $\Fp$ is mapped to $0$. Therefore, the map from $I_{(W,\psi)}(H\oplus\Fp,\eta\oplus\epsilon_{\Fp})$ to  $I_{(W,\psi)}(H,\eta)$ induced by the projection from $(H\oplus\Fp,\eta\oplus\epsilon_{\Fp})$ to $(H,\eta)$ is an isomorphism and $\decal I_{(W,\psi)}=0$.\\
    
    Since $I_V$ is analytic (cf \cite{K1}), $ \mathcal{U}(I_V)$ is analytic and therefore the tensor product $\IWV$ is analytic.
\end{proof}

Since the $\IWV$ form a family of injective cogenerators, any simple object $S$ embeds in some $\IWV$. Also, since $\IWV$ is analytic, $S$ embeds in some $p_n(\IWV)$ and is therefore polynomial.\\

An important feature of the category $\Sector$ is that, when there is a map $\gamma$ from $(H,\eta)$ to $(W,\psi)$ either $(H/\ker(\eta),\Tilde{\eta})$ and $(W/\ker(\psi),\Tilde{\psi})$ are isomorphic or there is no map from $(W,\psi)$ to $(H,\eta)$. Therefore, for $(W,\psi)$ a maximal element among isomorphism classes of objects in $\mathfrak{R}_S$ such that there exist $V$ with $F(W\oplus V,\psi\boxplus\epsilon_V)\neq 0$, one can consider the subfunctor $\Bar{F}$ of $F$, with $\Bar{F}(H,\eta)=F(H,\eta)$, if $(H/\ker(\eta),\Tilde{\eta})\cong (W,\psi)$ and $\Bar{F}=0$ otherwise. We get that, for $S$ simple, there is $(W,\psi)\in\mathfrak{R}_S$ such that $S(H,\eta)$ non trivial implies that $(H/\ker(\eta),\Tilde{\eta})\cong(W,\psi)$.\\

We can now describe the simple objects of $\Funcat$.\\

Let $S$ be a simple polynomial functor of degree $n$, $\decal^n$ maps $S$ onto a simple object of $\mathcal{F}(\mathfrak{R}_S,\Fp\left[\Sigma_n\right]-\mathcal{M}od)$. Those are the functors that map some $(W,\psi)\in\mathfrak{R}_S$ to some simple object in $\mathcal{F}(\mathcal{A}ut_{\Sector}(W,\psi),\Fp\left[\Sigma_n\right]-\mathcal{M}od)\cong\Fp\left[\mathcal{A}ut_{\Sector}(W,\psi)\times\Sigma_n\right]-\mathcal{M}od$, and the $(H,\eta)$ non isomorphic to $(W,\psi)$ to $0$.\\

In the following, we will use standard results about simple objects of $\mathcal{F}(\E,\Ef)$. We use the notations of \cite{PS}. For every $2$-regular partition $\lambda$, there is an element $\epsilon_\lambda\in\Fp\left[\Sigma_n\right]$, denoted $\Bar{R}_\lambda\Bar{C}_\lambda\Bar{R}_\lambda$ in \cite{PS}, such that $\epsilon_\lambda\Fp\left[\Sigma_n\right]$ is isomorphic to the simple module parametrized by $\lambda$.\\ 

It is known (cf \cite{PS}), that for $\epsilon_\lambda T^n$ the functor on $\E$ that maps $V$ to the image of $V^{\otimes n}$ under the right action of $\epsilon_\lambda$, $\epsilon_\lambda T^n$ is a polynomial functor of degree $n$ in $\mathcal{F}(\E,\Ef)$ and admits no non-trivial subfunctor of degree less than $n-1$. It is also known that $\Bar{\Delta}^n(\epsilon_\lambda T^n)\cong\epsilon_\lambda\Fp\left[\Sigma_n\right]$. 

\begin{theorem}\label{mainx}
    There is a one-to-one correspondence between isomorphism classes of simple objects of $\Funcat$ and isomorphism classes of simple objects of $$\bigsqcup\limits_{(W,\psi),n}\Fp\left[\mathcal{A}ut_{\Sector}(W,\psi)\times\Sigma_n\right]-\mathcal{M}od$$ with $(W,\psi)$ running through the isomorphism classes of objects in $\mathfrak{R}_S$ and $n$ running through $\N$.
\end{theorem}

\begin{proof}
We have already described the map from simple objects in $\Funcat$ to simple objects in $\bigsqcup\limits_{(W,\psi),n}\mathcal{F}(\mathcal{A}ut_{\Sector}(W,\psi),\Fp\left[\Sigma_n\right]-\mathcal{M}od)$. We have to prove that it is a one to one correspondence.\\

Let $M$ be a simple object in $\mathcal{F}(\mathcal{A}ut_{\Sector}(W,\psi),\Fp\left[\Sigma_n\right]-\mathcal{M}od)$. Since $\mathcal{A}ut_{\Sector}(W,\psi)$ is a category with only one object, $M$ is a $\Fp\left[\Sigma_n\right]$-module equipped with a left action of $\mathcal{A}ut_{\Sector}(W,\psi)$. As a $\Fp\left[\Sigma_n\right]$-module, it admits an injection from some simple $\Sigma_n$-module $\epsilon_\lambda\Fp\left[\Sigma_n\right]$. Each element of $\mathcal{A}ut_{\Sector}(W,\psi)$ maps $\epsilon_\lambda\Fp\left[\Sigma_n\right]$ to some isomorphic $\Fp\left[\Sigma_n\right]$-submodule of $M$. Those are either disjoint or equal to each-other. Using that $M$ is simple as an object in $\mathcal{F}(\mathcal{A}ut_{\Sector}(W,\psi),\Fp\left[\Sigma_n\right]-\mathcal{M}od)$, we get an isomorphism of $\Sigma_n$-modules $M\cong(\epsilon_\lambda\Fp\left[\Sigma_n\right])^{\oplus i}$ for some $i\in\N$ (this is because $\mathcal{A}ut_{\Sector}(W,\psi)$ is finite). We consider $T^n\otimes_{\Sigma_n}M\in\Funcat$. It admits a quotient of the form $(\epsilon_\lambda T^n)^{\oplus i}$ (by abuse of notation, we omit the action of morphisms in $\mathfrak{R}_S$ from the notation). This subfunctor admits no sub-functor of degree less than or equal to $n-1$ and $\decal^n(\epsilon_\lambda T^n)^{\oplus i}\cong M$, hence it is the quotient of $T^n\otimes_{\Sigma_n}M$ by $p_{n-1}(T^n\otimes_{\Sigma_n}M)$. Therefore, $(\epsilon_\lambda T^n)^{\oplus i}$ is a simple object in $\Funcat$.\\

Furthermore, let $F\in\Funcat$, polynomial of degree $n$ such that $\decal^n F\cong M$. The unit of the adjunction between $\decal^n$ and $T^n\otimes_{\Sigma_n}\_$, gives us a map from $T^n\otimes_{\Sigma_n}M$ to $F$, since it is an isomorphism in $\pol{n}/\pol{n-1}$, its kernel is included in $p_{n-1}(T^n\otimes_{\Sigma_n}M)$. Therefore, it factorises $T^n\otimes_{\Sigma_n}M\twoheadrightarrow(\epsilon_\lambda T^n)^{\oplus i}$. We get that either $F$ is not simple or it is isomorphic to $(\epsilon_\lambda T^n)^{\oplus i}$.
\end{proof}

\bibliographystyle{alpha}
\bibliography{Biblio}

\begin{thebibliography}{HLS93}

\bibitem[Blo23]{B2}
Ouriel Bloed\'e.
\newblock Approximation of the centre of unstable algebras using the nilpotent filtration, 2023.

\bibitem[Dja06]{djament:tel-00119635}
Aur{\'e}lien Djament.
\newblock {\em {Repr{\'e}sentations g{\'e}n{\'e}riques des groupes lin{\'e}aires : cat{\'e}gories de foncteurs en grassmanniennes, avec applications {\`a} la conjecture artinienne}}.
\newblock Theses, December 2006.

\bibitem[DV15]{DV}
Aur{\'e}lien Djament and Christine Vespa.
\newblock {Sur l'homologie des groupes d'automorphismes des groupes libres {\`a} coefficients polynomiaux}.
\newblock {\em {Commentarii Mathematici Helvetici}}, 90(1):33--58, 2015.

\bibitem[HLS93]{HLS2}
Hans-Werner Henn, Jean Lannes, and Lionel Schwartz.
\newblock The categories of unstable modules and unstable algebras over the {S}teenrod algebra modulo nilpotent objects.
\newblock {\em Math. Ann.}, 115:(1053--1106), 1993.

\bibitem[Kuh94]{K1}
Nicholas~J. Kuhn.
\newblock Generic representations of the finite general linear group and the {S}teenrod algebra : {I}.
\newblock {\em American Journal of Mathematics}, 116:(327--360), 1994.

\bibitem[Pir88]{P1}
Teimuraz Pirashvili.
\newblock Polynomial functors.
\newblock {\em Trudy Tbiliss. Mat. Inst. Razmadze Akad. Nauk Gruzin}, 91:55 -- 66, 1988.

\bibitem[PS98]{PS}
Laurent Piriou and Lionel Schwartz.
\newblock Extensions de foncteurs simples.
\newblock {\em K-Theory}, 15, 1998.

\end{thebibliography}

\begin{center}
Address : Université de Lille, laboratoire Paul Painlevé (UMR8524), Cité scientifique, bât. M3, 59655
VILLENEUVE D’ASCQ CEDEX, FRANCE\\
    e-mail : aacde13@live.fr
    \end{center}

\end{document}